\magnification=\magstep1
\input amstex
\documentstyle{amsppt}

\define\defeq{\overset{\text{def}}\to=}
\define\ab{\operatorname{ab}}
\define\pr{\operatorname{pr}}
\define\Gal{\operatorname{Gal}}

\define\sep{\operatorname{sep}}
\def \isom {\overset \sim \to \rightarrow}

\define\Br{\operatorname{Br}}

\define\id{\operatorname{id}}

\def \Br{\operatorname{Br}}
\def \f{\operatorname{f}}
\def \Sel{\operatorname{Sel}}
\def \Im{\operatorname{Im}}

\NoRunningHeads
\NoBlackBoxes
\topmatter

\title
On the  birational anabelian section conjecture
\endtitle

\author
Mohamed Sa\"\i di
\endauthor

\abstract Assuming the finiteness of the Shafarevich-Tate group of elliptic curves over number fields we make several observations
on the birational Grotendieck anabelian setion conjecture.
We prove that the birational setion conjecture for curves over number fields can be reduced to the case of elliptic curves.
In this case we prove that, as a consequence of a result of Stoll, a section of the exact sequence of the absolute Galois group of an elliptic curve over a number field arises from a rational point
if and only if the induced section of the corresponding (geometrically abelianised) arithmetic fundamental group of the elliptic curve arises from a rational point. 
We also prove that given any curve over a number field, there exists a double covering of this curve for which the birational setion conjecture holds true.
\endabstract
\toc

\subhead
\S 0. Introduction
\endsubhead

\subhead
\S 1. The Birational Grothendieck Anabelian Section Conjecture
\endsubhead

\subhead
\S 2. Reduction of the Birational Section Conjecture to Curves with a given Genus
\endsubhead

\subhead
\S 3. Birational Sections for Genus $0$ Curves over Number Fields
\endsubhead

\subhead
\S 4. Birational Sections for Genus $1$ Curves over Number Fields
\endsubhead

\subhead
\S 5. Birational Sections for Genus $g\ge 2$ Curves over Number Fields
\endsubhead

\endtoc

\endtopmatter
\document

\subhead
\S 0. Introduction
\endsubhead
Let $k$ be a field of characteristic $0$, and $X$ a proper, smooth, and geometrically connected
(not necessarily hyperbolic) {\it algebraic curve} over $k$. Let $K_X$ be the function field of $X$, 
$K_X^{\sep}$ a separable closure of $K_X$, and $\bar k$ the algebraic closure of $k$ in
$K_X^{\sep}$.  
Write 
$$G_X\defeq \Gal (K_X^{\sep}/K_X),$$
and 
$$G_{\overline X}\defeq  \Gal (K_X^{\sep}/K_{\overline X}),$$
where $K_{\overline X}\defeq K_X.\overline k$ is the function field of the
geometric fibre $\overline X\defeq X\times _k\overline k$ of $X$.
There exists a canonical exact sequence of profinite absolute Galois groups
$$1\to G_{\overline X}\to G_X@>\pr_X>> G_k\to 1,$$
where $G_k=\Gal (\bar k/k)$.
Let $x\in X(k)$ be a {\it rational} point. Then $x$ determines a {\it decomposition subgroup}
$D_x\subset G_X$, which is only defined up to conjugation by the elements of $G_{\overline X}$, and which maps surjectively onto 
$G_k$ via the natural projection $\pr_X:G_X\twoheadrightarrow G_k$. 
More precisely, $D_x$ sits naturally in the following exact sequence
$$1\to \hat \Bbb Z(1)\to D_x\to G_k\to 1.$$
The above exact sequence is known to be split. A {\it section} $G_k\to D_x$ of the natural projection $D_x\twoheadrightarrow G_k$ (i.e. a splitting of the above exact sequence)
determines naturally a {\it section} $G_k\to G_X$ of the natural projection $\pr _X:G_X\twoheadrightarrow G_k$, whose image is contained in $D_x$.
The {\it birational} version of the anabelian Grothendieck section conjecture for curves predicts that a {\it section}, or {\it splitting}, of the exact sequence 
$$1\to G_{\overline X}\to G_X@>\pr_X>> G_k\to 1,$$
over a {\it finitely generated} field $k$ of characteristic zero, necessarily arises from a {\it rational point} $x\in X(k)$ of the curve $X$ as explained above. (cf. [Koenigsmann], and $\S1$, for more details).  
More generally, one says that a field $k$ possesses the {\it birational section property for curves} if a similar statement as above holds for any curve over $k$ (cf. loc. cit.).
A major breakthrough towards the birational section conjecture is the fundamental result of Koenigsmann, that $p$-adic local fields (i.e. finite extensions of $\Bbb Q_p$)
possess the birational section property for curves (cf. [Koenigsmann]).  Also, it is well-known that the field $\Bbb R$ of real numbers has the birational section property.

In this paper we make several observations, and prove several facts,  regarding this conjecture. First, we prove that 
in order to verify  that a field $k$ has the birational section property it suffices to reduce to the case where $X=\Bbb P^1_k$ is the projective line 
(cf. Lemma 2.1), and more generally to the case of curves with a given genus $g\ge 1$ 
(cf. Proposition 2.2 and Corollary 2.3).  

Let $k$ be a number field. For each place $v$ of $k$ let $k_v$ be the completion of $k$ at $v$
and $X_v\defeq X\times_kk_v$. Let $s:G_k\to G_X$ be a section of the natural projection $G_X\twoheadrightarrow G_k$. Then 
$s$ gives {\it naturally} rise to sections $s_v:G_{k_v}\to G_{X_v}$ of the natural projection $G_{X_v}\twoheadrightarrow G_{k_v}$,
for each place $v$ of $k$ (cf. Proof of Proposition 1.4).  Each of these sections $s_v:G_{k_v}\to G_{X_v}$ arises from a rational point $x_v\in X(k_v)$,
since $p$-adic local fields and the field of real numbers  possess  the birational section property for curves. We prove the following.
Suppose there exists a rational point $x\in X(k)$ such that $x=x_v$ for each place $v$ of $k$, i.e. the local sections $s_v$ arise form a global rational point $x\in X(k)$,
then the section $s:G_k\to G_X$ arises from the rational point $x$ (cf. Proposition 5.3, and the Proof of Proposition 3.3 where a more precise 
statement is proved in the case where $X=\Bbb P^1_k$ is the projective line).

Assuming the {\it finiteness} of the Shafarevich-Tate groups of elliptic curves, we prove that the birational section conjecture for curves over number fields can be reduced to the case of {\it elliptic curves} (cf. Proposition 4.2). In the case of an {\it elliptic curve} $E$ over a {\it number field} $k$, with {\it finite Shafarevich-Tate group}, a section
$s:G_k\to G_E$
 of the  exact sequence 
$$1\to G_{\overline E}\to G_E@>\pr_E>> G_k\to 1$$
of the absolute Galois group of the function field of $E$ gives rise naturally to a section $\tilde s: G_k\to \Pi_E$ of the exact sequence 
$$1\to T\overline E \to \Pi_E@>\pr_E>> G_k\to 1$$
of the arithmetic fundamental group $\Pi_E$ of $E$, where $T\overline E$ is the Tate module of $E$.
Fix a base point of the torsor of splittings of the exact sequence $1\to T\overline E \to \Pi_E@>\pr_E>> G_k\to 1$ which arises from the origin of $E$.
then the {\it conjugacy class} of the section $s:G_k\to \Pi_E$ corresponds to an element of $H^1(G_k,T\overline E)$ which we denote also $\tilde s$.
We observe that $\tilde s$ lies in the subgroup (via Kummer theory) $E(k)^{\wedge}$ of $H^1(G_k,T\overline E)$, where $E(k)^{\wedge}$ denotes the profinite completion of the 
group of rational points $E(k)$ (cf. Lemma 4.4). Furthermore, we prove that, as a consequence of a result of Stoll, the above birational section $s:G_k\to G_E$ arises from
a rational point of $E$ if and only the above element $\tilde s \in E(k)^{\wedge}$ lies in the discrete subgroup $E(k)\subset E(k)^{\wedge}$ (cf. Propsoition 4.6). More precisely, one can in the framework of the birational anabelian section conjecture give a (group-theoretic) characterisation of the {\it discrete} group $E(k)$ inside its {\it profinite} completion 
$E(k)^{\wedge}$ (cf. Proposition 4.7).
Similar observations are made for birational sections in the case of curves of genus at least $2$ (cf. Propsoition 5.2).
Finally, we prove that given a proper, smooth, and geometrically connected curve $X$ over a number field $k$ there exists
a {\it double covering} $X'\to X$ defined over $k$ such that the birational section conjecture holds true for $X'$, under the assumption that the
Shafarevich-Tate groups of elliptic curves over $k$ are finite, (cf. Lemma 5.5, and Remark 5.6).

\definition{Acknowledgment} I would like to thank Akio Tamagawa for the discussions we had  around the topic of this paper, and for 
his comments on an earlier version of this paper. I also thank him for inviting me to the 
Research Institute for Mathematical sciences RIMS of Kyoto university where part of this work was done.
I also thank Alexei Skorobogatov for the discussion we had around the topic of this paper.

 \enddefinition

\subhead
\S 1. The Birational Grothendieck Anabelian Section Conjecture
\endsubhead
In this section we briefly recall, and explain, the content of the birational anabelian section conjecture of Grothendieck for curves (cf. [Grothendieck]).
We also fix notations that will be used throughout this paper.

Let $k$ be a field of characteristic $0$, and $X$ a proper, smooth, and geometrically connected
(not necessarily hyperbolic) {\it algebraic curve} over $k$. Let $K_X$ be the function field of $X$, 
$K_X^{\sep}$ a separable closure of $K_X$, and $\bar k$ the algebraic closure of $k$ in
$K_X^{\sep}$.  
Write 
$$G_X\defeq \Gal (K_X^{\sep}/K_X),$$
and 
$$G_{\overline X}\defeq  \Gal (K_X^{\sep}/K_{\overline X}),$$
where $K_{\overline X}\defeq K_X.\overline k$ is the function field of the
geometric fibre $\overline X\defeq X\times _k\overline k$ of $X$.
There exists a canonical exact sequence of profinite absolute Galois groups

$$1\to G_{\overline X}\to G_X@>\pr_X>> G_k\to 1, \tag {$1$}$$
where $G_k=\Gal (\bar k/k)$.

By a group-theoretic {\it section}, or a {\it splitting}, of the exact sequence (1) we mean a continuous homomorphism
$s:G_k\to G_X$ such that $\pr_X\circ s=\id_{G_k}$.
Let $x\in X(k)$ be a {\it rational} point of $X$.  Then $x$ determines a {\it decomposition subgroup}
$D_x\subset G_X$, which is only defined up to conjugation by the elements of $G_{\overline X}$, and which maps surjectively onto 
$G_k$ via the natural projection $\pr_X:G_X\twoheadrightarrow G_k$. 
More precisely, $D_x$ sits naturally in the following exact sequence
$$1\to \hat \Bbb Z(1)\to D_x\to G_k\to 1.\tag {$2$}$$
The exact sequence (2) is known to be split. Indeed, the extension defined by extracting $n$-th roots, for all positive integers $n$, of a given local parameter at $x$
defines a splitting of this sequence. The set of all splittings of the exact sequence (2) is a torsor under the Galois cohomology group $H^1(G_k,\hat \Bbb Z(1))$. 
A section $G_k\to D_x$ of the natural projection $D_x\twoheadrightarrow G_k$ (i.e. a splitting of the exact sequence (2))
determines naturally a section $G_k\to G_X$ of the natural projection $\pr _X:G_X\twoheadrightarrow G_k$, whose image is contained in $D_x$.

\definition {The Birational Grothendieck Anabelian Section Conjecture (BGASC) (cf. [Koenigsmann])} 
Assume that $k$ is {\it finitely generated} over the
prime field $\Bbb Q$. Let $s:G_k\to G_X$ be a group-theoretic {\it section} of the natural projection
$\pr_X:G_X \twoheadrightarrow G_k$. Then the image $s(G_k)$ is contained in a decomposition subgroup $D_x\subset G_X$ associated to a {\it unique
rational point} $x\in X(k)$. In particular, the existence of the section $s$ implies that $X(k)\neq \varnothing$.
\enddefinition

\definition {Definition 1.1} Let $k$ be a field.  We say that the BGASC {\it holds true} over $k$ if for every proper, smooth, and geometrically connected
algebraic curve $X$ over $k$, and every group-theoretic {\it section}
$s:G_k\to G_X$ of the natural projection
$\pr_X: G_X \twoheadrightarrow G_k$, the image $s(G_k)$ is contained in a decomposition subgroup $D_x\subset G_X$ associated to a
{\it unique rational point} $x\in X(k)$. In this case we say that the BGASC {\it holds true for the $k$-curve} $X$.
\enddefinition

In connection with the BGASC, in the case where $k$ is a {\it number field}, it is natural to formulate a $p$-adic version of this conjecture over $p$-adic local fields.

\definition {A $p$-adic Version of the Birational Grothendieck Anabelian Section Conjecture ($p$-adic BGASC) (cf. loc. cit.)} 
Let $p> 0$ be a prime integer, and assume that $k$ is a {\it finite extension} of $\Bbb Q_p$. Then the BGASC holds true over $k$. More precisely,
let $s:G_k\to G_X$ be a group-theoretic {\it section} of the natural projection
$G_X \twoheadrightarrow G_k$. Then the image $s(G_k)$ is contained in a decomposition subgroup $D_x\subset G_X$ associated to a {\it unique
rational point} $x\in X(k)$. In particular, the existence of the section $s$ implies that $X(k)\neq \varnothing$.
\enddefinition

\definition {Remark 1.2} The uniqueness of the rational point $x\in X(k)$ mentioned in the BGASC, and its $p$-adic variant, is well-known if such a point exists. 
Indeed, any conjugates of two decomposition subgroups of $G_X$ corresponding to distinct closed points of $X$ have trivial intersection.
Thus,  in order to establish these conjectures, it suffices to establish the existence of a rational point $x$
such that a corresponding decomposition group $D_x$ contains the image of the section $s$.
\enddefinition

A major breakthrough towards the BGASC is the following fundamental result concerning the $p$-adic BGASC, and which is du to Koenigsmann (cf. [Koenigsmann]).

\proclaim {Theorem 1.3 (Koenigsmann)} The $p$-adic version of the BGASC holds true. More precisely, 
assume that $k$ is a finite extension of $\Bbb Q_p$.
Let $s:G_k\to G_X$ be a group-theoretic section of the natural projection
$G_X \twoheadrightarrow G_k$. Then the image $s(G_k)$ is contained in a decomposition subgroup $D_x$ associated to a
unique rational point $x\in X(k)$. In particular, the existence of the section $s$ implies that  $X(k)\neq \varnothing$.
\endproclaim

This result has been strengthened by Pop, who proved a $\Bbb Z/p\Bbb Z$-{\it meta-abelian} version of this theorem  (see [Pop] for more details). 
An important consequence of Theorem 1.3 is the following, which was already observed in [Koenigsmann].

\proclaim {Proposition 1.4} Let $k$ be a number field and $X$ a proper, smooth, and geometrically connected (not necessarily hyperbolic) curve over $k$.
Assume that there exists a section $s:G_k\to G_X$ of the natural projection $G_X\twoheadrightarrow G_k$.
Then the section $s$ gives rise to an adelic point  $(x_v)_v\in X(\Bbb A_k)$. 
Moreover, $x_v$ (resp. the connected component containing $x_v$) is uniquely determined by the section $s$ in the case where $v$ is a finite place (resp. $v$ is a real place).
Here,  $\Bbb A_k$ denotes the ring of ad\`eles of $k$, and $v$ runs over all places of $k$.
\endproclaim

\demo {Proof}
See [Koenigsmann], Corollary 2.6. In fact one can prove a more precise statement than in [Koenigsmann] (cf. loc. cit.).
For each place $v$ of $k$, let $k_v^h$ (resp. $k_v$) be the henselisation of $k$ at $v$ (resp. the completion of 
$k$ at $v$), and $X_v^h\defeq X\times _k k_v^h$ (resp. $X_v\defeq X\times _k k_v$). The section $s$ induces naturally a section $s_v^h:G_{k_v}\to G_{X_{v}^h}$
of the natural projection $G_{X_{v}^h}\twoheadrightarrow G_{k_v}$ (here, we fix an identification of $G_{k_v}\isom G_{k_v}^h$ with a decomposition subgroup of $G_k$ at the place $v$).
By a (an unpublished) result of Tamagawa, the section $s_v^h$ can be lifted to a section $s_v:G_{k_v}\to G_{X_{v}}$
of the natural projection $G_{X_{v}}\twoheadrightarrow G_{k_v}$ (cf. [Sa\"\i di], Theorem 5.6).  More precisely, one can construct a section $s_v:G_{k_v}\to G_{X_v}$ 
which fits into the following commutative diagram
$$
\CD
G_{k_v} @>s_v>>  G_{X_{v}} \\
@V{\id}VV      @VVV    \\
G_{k_v} @>s_v^h>>  G_{X_{v}^h}\\
@VVV    @VVV \\
G_k  @>s>> G_X \\
\endCD
$$
where the right top vertical map is a natural surjection, and the left low vertical map is an embeeding. 
By Theorem 1.3 above of Koenigsmann, for a finite place $v$, the image $s_v(G_{k_v})$ is contained in a decomposition subgroup
$D_{x_v}$ associated to a unique rational point $x_v\in X(k_v)$. This is also true for the archimedian places (the so-called real section conjecture holds true).
Note that in the case where $v$ is a real place only the connected component of $X(k_v)$ containing $x_v$ is well determined by the section $s$.
Thus, to the section $s$ is associated {\it naturally} an adelic point $(x_v)_v\in X(\Bbb A_k)$ with the required properties. 
\qed
\enddemo

\subhead
\S 2. Reduction of the Birational Section Conjecture to Curves with a given Genus
\endsubhead
In this section we state and prove our main observation concerning the BGASC, that it can be reduced to curves with a given genus. 
Our first observation is that the BGASC can be (easily) reduced to the case of the {\it projective line}.

\proclaim {Lemma 2.1}  Let $k$ be a field. Assume that the BGASC holds true for $\Bbb P^1_k$ (cf. Definition 1.1). Then the BGASC
holds true for any $k$-curve $X$ which is projective, smooth, and geometrically connected.
\endproclaim

\demo {Proof} Let $X$ be a projective, smooth,
and geometrically connected algebraic curve over $k$. Let $f:X\to \Bbb P^1_k$ be a finite morphism, which corresponds to a finite field extension
$K_X/k(T)$ where $k(T)\defeq K_{\Bbb P^1_k}$. We have a natural commutative diagram
of exact sequences of profinite Galois groups

$$
\CD
1 @>>> G_{\overline X} @>>>  G_X@>\pr_X>> G_k  @>>>  1 \\
@.        @VVV              @VVV           @V{\id}VV    \\
1@>>> G_{\Bbb P^1_{\overline k}}@>>> G_{\Bbb P^1_k} @>{\pr_{\Bbb P^1_k}} >> G_k  @>>>  1\\
\endCD
$$
where the left and middle  vertical maps are natural inclusions. Here, $G_{\Bbb P^1_k}\defeq \Gal (K_X^{\sep}/k(T))$,
and $G_{\Bbb P^1_{\overline k}}\defeq \Gal (K_X^{\sep}/\overline k(T))$.
Let $s:G_k\to G_X$ be a group-theoretic section of the 
natural projection $\pr_X:G_X\twoheadrightarrow G_k$. The image of $s(G_k)$ in $G_{\Bbb P^1_k}$, 
via the natural embedding $G_X\hookrightarrow G_{\Bbb P^1_k}$,  determines a group-theoretic section $\tilde s:G_k\to
G_{\Bbb P^1_k}$ of the natural projection $G_{\Bbb P^1_k} \twoheadrightarrow G_k$. Assume that the BGASC holds
true for $\Bbb P^1_k$. Then the image $\tilde s(G_k)$ of the section $\tilde s$ is contained in a decomposition subgroup
$D_y\subset G_{\Bbb P^1_k}$ associated to a unique rational point $y\in \Bbb P^1(k)$. The intersection $D_x\defeq D_y \cap G_X$ is then
the decomposition group associated to a unique point $x\in X$, which is necessarily $k$-rational since $D_x$
maps surjectively onto $G_k$ via the natural projection $G_X\twoheadrightarrow G_k$.    Moreover, we have $s(G_k)\subseteq D_x$.
\qed
\enddemo

Next, we prove that the BGASC over any field can be reduced to the case of curves with a given genus $g\ge 1$.
We refer to the discussion in $\S 3$ for the case of genus $0$ curves over number fieds.

\proclaim {Proposition 2.2}  Let $k$ be a field, and $g\ge 1$ an integer. 
Assume that the BGASC holds true for genus $g$ proper, smooth, and geometrically connected curves over $k$. 
Then the BGASC holds true for the projective line over $k$. 
\endproclaim

In particular, as a consequence of Lemma 2.1, and Proposition 2.2,  one deduces immediately the following.

\proclaim {Corollary 2.3}  Let $k$ be a finitely generated field over $\Bbb Q$, and $g\ge 1$ an integer. 
Assume that the
BGASC holds true for genus $g$  proper, smooth, and geometrically connected curves over $k$. Then the BGASC holds true for any projective, smooth, and geometrically connected curve $X$ over $k$. 
\endproclaim

\demo {Proof of Proposition 2.2} 
Recall the exact sequence of absolute Galois groups
$$1\to  G_{\Bbb P^1_{\overline k}}\to  G_{\Bbb P^1_k} @>{\pr_{\Bbb P^1_k}} >> G_k  \to  1.$$
Let $s:G_k\to G_{\Bbb P^1_k}$ be a section of the natural projection $G_{\Bbb P^1_k}\twoheadrightarrow G_k$.
Let $\overline \Delta$ be an open subgroup of $G_{\Bbb P^1_{\overline k}}$ corresponding to a finite morphism 
$\tilde f:\Tilde X\to \Bbb P^1_{\overline k}$, where $\Tilde X$ is a genus $g$ proper, smooth, and connected curve over $\overline k$.  
Assume moreover that the finite morphism $\tilde f:\Tilde X\to \Bbb P^1_{\overline k}$ is defined over $k$, in which case $\overline \Delta$ is stable under the natural action of $s(G_k)$ on $G_{\Bbb P^1_{\overline k}}$ via inner automorphisms.
Write $\Delta\defeq \overline \Delta . s(G_k)$. Then $\Delta$ is an open subgroup of
$G_{\Bbb P^1_{k}}$ which corresponds to a finite morphism $f:X\to \Bbb P^1_k$ where $X$ is a projective, smooth, and geometrically connected
$k$-curve. Let $G_X\defeq \Gal (K^{\sep}/K_X)=\Delta $, and $G_{\overline X}\defeq  \Gal (K^{\sep}/K_{\overline X})=\overline \Delta$, where $\overline X\defeq X\times _k
\overline k$, and $K^{\sep}\defeq K_{\Bbb P^1_k}^{\sep}$. We have a natural commutative diagram of exact sequences of absolute Galois groups

$$
\CD
1 @>>> G_{\overline X}=\overline \Delta @>>>  G_X=\Delta @>\pr_X>> G_k  @>>>  1 \\
@.        @VVV              @VVV           @V{\id}VV    \\
1@>>> G_{\Bbb P^1_{\overline k}}@>>> G_{\Bbb P^1_k} @>{\pr_{\Bbb P^1_k}} >> G_k  @>>>  1\\
\endCD
$$
Note that by construction we have a natural isomorphism $\Tilde X\isom \overline X$ over $\overline k$.
In particular, $X$ is a genus $g$ curve.
Also, by construction, the group-theoretic section $s:G_k\to G_{\Bbb P^1_k}$ naturally restricts to a group-theoretic section $s:G_k\to G_X$ of the natural projection
$G_X\twoheadrightarrow G_k$.  Furthermore, in order to show that  $s(G_k)$ is contained in a decomposition subgroup $D_y$ associated to a rational point 
$y\in \Bbb P^1(k)$, it suffices to show that $s(G_k)$ is contained in a decomposition subgroup $D_x\subset G_X$ associated to a rational point $x\in X(k)$. Indeed, if 
$s(G_k)\subseteq D_x$, where $x\in X(k)$, then $s(G_k)\subseteq D_y$ is contained in a decomposition subgroup associated to the rational point 
$y\in \Bbb P^1(k)$ which is the image of $x$ under the morphism $f:X\to \Bbb P^1_k$ corresponding to the inclusion $G_X\subset G_{\Bbb P^1_k}$. 
Moreover, $s(G_k)\subset D_x$ for a unique rational point $x\in X(k)$ if we assume that the BGASC holds true for $X$. 
\qed
\enddemo

\subhead
\S 3. Birational Sections for Genus $0$ Curves over Number Fields
\endsubhead
In this section we discuss the BGASC in the case of genus $0$ curves, mainly over number fields. We first observe the following.

\proclaim {Proposition 3.1} Let $k$ be a number field and $X$ a proper, smooth, and geometrically connected genus $0$ curve over $k$.
Assume that there exists a section $s:G_k\to G_X$ of the natural projection $G_X\twoheadrightarrow G_k$. Then $X(k)\neq \varnothing$. In particular,
$X\isom \Bbb P^1_k$ is a projective line.   
\endproclaim

\demo {Proof}
Indeed, the set of adelic point $X(\Bbb A_k)\neq \varnothing$ is non-empty by Proposition 1.4. Hence the set of rational points
$X(k)\neq \varnothing$ is non empty, since the Hasse 
principle for rational points holds for $X$.
\qed
\enddemo

\definition {3.2} Let $k$ be a number field, and $X=\Bbb P^1_k$ the projective line over $k$. Let $\infty \in X(k)=\Bbb P^1_k(k)$ be a rational point.
Let $s:G_k\to G_X$ be a section of the natural projection $G_X\twoheadrightarrow G_k$. 
For each place $v$ of $k$, let $k_v$ be the completion of 
$k$ at $v$, and $X_v\defeq X\times _k k_v$. The section $s$ induces naturally a section $s_v:G_{k_v}\to G_{X_{v}}$
of the natural projection $G_{X_{v}}\twoheadrightarrow G_{k_v}$ (here, we fix an identification of $G_{k_v}$ with a decomposition subgroup of $G_k$ at the place $v$).
More precisely, there exists a section $s_v$ which fits into the following 
commutative diagram
$$
\CD
G_{k_v} @>s_v>>  G_{X_{v}} \\
@VVV      @VVV    \\
G_k @>s_v^h>>  G_{X} \\
\endCD
$$
where the right vertical map is the natural one (cf. Proof of Proposition 1.4). 
We know that the section $s_v$ arises from a unique rational point $x_v\in X(k_v)$ (cf. Theorem 1.3).
If the section $s$ arises from a rational point $x\in X(k)$. Then, after observing the natural action of $PGL_2(k)$ on $G_X$, we can assume that $x=\infty$.
The section $s_v$ would then also arise from the point $\infty$.
Reciprocally, We can prove the following.
\enddefinition

\proclaim {Proposition 3.3} We use the same notations and assumptions as in 3.2. Assume that for each place $v$ of $k$ we have $x_v=\infty$. In other words assume 
the image $s_v(G_{k_v})\subset D_{\infty,v}$ is contained in a decomposition group $D_{\infty,v}\subset G_{X_v}$ associated to the point $\infty \in X(k)\subset X(k_v)$.
Then the section $s$ arises from the rational point $\infty$, i.e. the image $s(G_k)\subset D_{\infty}$ is contained in a decomposition group $D_{\infty}\subset G_k$
associated to $\infty$.
\endproclaim

\demo {Proof} 
One can write down a proof similar to the proof of Proposition 4.6, which resorts directly to Theorem 1.3 and a result of Stoll (cf. the Proof of Proposition 4.6).
We will however prove a slightly more precise statement without resorting to the result of Stoll.
We will show that there exists a {\it neighbourhood} of the section $s$, that is  the absolute Galois group of
the function field of an elliptic curve for which the BGASC holds (under the assumption of finiteness of Shafarevich-Tate groups of elliptic curves over $k$).

Let $E$ be an elliptic curve over $k$ with trivial Mordell-Weil rank (such curves exist, cf. Remark 5.6). 
Then $E$ can be realised as a Galois cover $f:E\to \Bbb P^1_k$ of degree $2$
of the projective line ramified above $\infty$ (in particular, the rational point $\infty$ lifts to a unique rational point $x\in E(k)$), and the absolute Galois group $G_E$ 
naturally embeds $\iota :G_E\hookrightarrow G_{\Bbb P^1_k}$ in $G_{\Bbb P^1_k}$ as a normal subgroup of index $2$. 
We have a natural commutative diagram
$$
\CD
G_{E} @>>>  G_{k} \\
@V{\iota}VV      @V{\id}VV    \\
G_{\Bbb P^1_k} @>>>  G_{k} \\
\endCD
$$
We will show that $G_E$ necessarily contains the image $s(G_k)$ of the section $s$. For each place $v$ of $k$ we have a natural commutative diagram

$$
\CD
G_{E_v} @>>>  G_{k_v} \\
@V{\iota_v}VV      @V{\id}VV    \\
G_{\Bbb P^1_{k_v}}@>>>  G_{k_v} \\
\endCD
$$
where the left vertical embedding is naturally induced by $\iota$. Moreover, the image $s_v(G_{k_v})$ of the section $s_v$, and all its conjugate, 
are contained in $G_{E_v}$ (since
the point $\infty$ lifts to a unique rational point of $E_v$,
and $G_{E_v}$ is a normal subgroup of $G_{\Bbb P^1_{k_v}}$). On the other hand, $G_k$ is normally topologically generated by the decomposition subgroups $G_{k_v}$ as follows from the Chebotarev density theorem. From this follows that $G_E$ is normally topologically generated by the images of the $G_{E_v}$, where $v$ runs over all places of $k$. Hence, $G_E$ contains $s(G_{k_v})$, and all its conjugates, for all places $v$. Thus, $G_E$ contains $s(G_k)$, and the section $s$ naturally restricts to a section
$s:G_k\to G_E$ of the natural projection $G_E\twoheadrightarrow G_k$, which arises from a rational point $y\in E(k)$ by Corollary 4.8 (here we assume that the Shafarevich-Tate group of $E$ is finite). In particular, the section $s:G_k\to G_{\Bbb P^1_k}$ arises from the rational point $x\in \Bbb P^1_k(k)$ which is the image of $y$ under the above morphism $f:E\to \Bbb P^1_k$. 
\qed
\enddemo

\definition {Remark/Question 3.4}
We use the same notations as in 3.2.
One can, after observing the action of $PGL_2(k)$, assume that for every finite set of places $S$ of $k$ one has $x_v=\infty$, 
since $PGL_2(k)$ is dense in $PGL_2(\Bbb A_k)$. Is it possible to prove that this leads to the same conclusion as in Proposition 3.3?
If yes, this would prove the BGASC for $\Bbb P^1_k$ in the case where $k$ is a number field.
\enddefinition

\subhead
\S 4. Birational Sections for Genus $1$ Curves over Number Fields
\endsubhead
In this section we discuss the BGASC for genus $1$ curves, mainly over number fields.
First, we observe that the existence of a birational section for a genus $1$ curve over a number field implies, assuming the 
finiteness of the Shafarevich-Tate groups for elliptic curves, that this curve is an elliptic curve. More precisley, we have the following.

\proclaim {Proposition 4.1} Let $k$ be a number field and $X$ a proper, smooth, geometrically connected genus $1$ curve over $k$. 
Assume that the Shafarevich-Tate groups of elliptic curves over $k$ are finite. Assume there exists a section $s:G_k\to G_X$
of the natural projection $G_X\twoheadrightarrow G_k$. Then $X(k)\neq 0$. In particular, $X$ is an elliptic curve. 
\endproclaim

An immediate consequence of Proposition 4.1, and Proposition 2.2, is that
the BGASC for {\it curves over number fields} can be reduced to the case of {\it elliptic curves over number fields} (assuming
the finiteness of the Shafarevich-Tate groups for elliptic curves). More precisely, we have the following. 

\proclaim {Proposition 4.2}  Let $k$ be a number field. 
Assume that the Shafarevich-Tate groups of elliptic curves over $k$ are finite, and that the
BGASC holds true for all elliptic curves over $k$. Then the BGASC holds true for any projective, smooth, and geometrically connected curve $X$ over $k$. 
\endproclaim

\demo {Proof of Proposition 4.1} Recall the exact sequence of absolute Galois groups
$$1\to  G_{\overline X}\to  G_{X} @>{\pr_{X}} >> G_k  \to  1.$$
Let $s:G_k\to G_{X}$ be a section of the natural projection $G_{X}\twoheadrightarrow G_k$.
By assumption, $X$ is a genus $1$ curve. Moreover, $X$ is a principal homogeneous space
over $k$ under its jacobian $E'$ which is an elliptic curve over $k$, and corresponds to an element of the Galois cohomology group $H^1(G_k,E')$. 
Next, assuming that the Shafarevich-Tate group of $E'$ is finite, we will show that
$X\isom E'$ is an elliptic curve. The existence of the section $s:G_k\to G_X$ implies that $X(\Bbb A_k)\neq \varnothing$ (cf. Proposition 1.4).
Thus, to the section $s$ is associated an adelic point $(x_v)_v\in X(\Bbb A_k)$ (cf. loc. cit.). 
The adelic point $(x_v)_v\in X(\Bbb A_k)$ survives every finite \'etale abelian descent obstruction
(cf. [Stoll], Definition 5.2), as follows easily from the existence of the global section $s$ (see also [Harari-Stix], Proposition 1.1), i.e. $(x_v)_v\in X(\Bbb A_k)^{\f-\ab}$
in the terminology of Stoll, where $X(\Bbb A_k)^{\f-\ab}$ is the set of adelic points cut out by the finite \'etale abelian descent conditions (cf. loc. cit. Definition 5.4).
For a different argument to deduce the existence of a point in $X(\Bbb A_k)^{\f-\ab}$ one may also use similar arguments as in the proof of 
Theorem 3.2 in [Harari-Stix]. On the other hand one has the following equality $X(\Bbb A_k)^{\f-\ab}=X(\Bbb A_k)^{\Br}$ where $X(\Bbb A_k)^{\Br}$
denotes the set of adelic points cut out by the Brauer-Manin conditions, i.e. the Brauer-Manin set (cf. [Stoll], Corollary 7.3). The non-emptiness of the Brauer set
$X(\Bbb A_k)^{\Br}$ implies, under the assumption that the Shafarevich-Tate group of $E'$ is finite, that $X(k)\neq \varnothing$ by a result of Manin (cf. [Manin]), 
hence $X\isom E'$ is an elliptic curve.
\qed
\enddemo

\definition {4.3} Next, we will discuss the BGASC in the case of an elliptic curve over a number field. In what follows we will assume that $k$ is a {\bf number field}, and $E$ is an {\bf elliptic curve over $k$ with finite Shafarevich-Tate group}. 
\enddefinition
Recall the exact sequence
of absolute Galois groups
$$1\to  G_{\overline E}\to  G_{E} @>{\pr_{E}} >> G_k  \to  1.$$
Let $\Pi_E$ be the quotient of $G_E$ which corresponds to the maximal {\it everywhere unramified} extension of $K_E$ contained in $K_E^{\sep}$.
Thus, $\Pi_E$ is the {\it arithmetic \'etale fundamental group} of $E$.
We have a natural commutative diagram of exact sequences
$$
\CD
1 @>>> G_{\overline E} @>>>  G_E@>\pr_E>> G_k  @>>>  1 \\
@.        @VVV              @VVV           @V{\id}VV    \\
1@>>> \Pi_{\overline E}@>>> \Pi_E @>{\pr_E} >> G_k  @>>>  1\\
\endCD
$$
where $\Pi_{\overline E}$ is the \'etale fundamental group of $\overline E\defeq E\times _k\bar k$, which is naturally identified with the {\it Tate module} $T\overline E$ of $\overline E$.
The left and middle vertical maps in the above diagram are surjective.
We fix a {\it base point} of the torsor of splittings of the exact sequence 
$1\to \Pi_{\overline E}\to  \Pi_E @>{\pr_E} >> G_k  \to 1$, which corresponds to the splitting arising from the origin of the elliptic curve $E$. 
The set of splittings of the above sequence is then a torsor under the Galois cohomology group $H^1(G_k,\Pi_{\overline E})$.
Let $s:G_k\to G_E$ be a group-theoretic {\it section} of the natural projection $\pr_E:G_E\twoheadrightarrow G_k$. Then $s$ induces naturally a group-theoretic section
$\tilde s:G_k\to \Pi_E$ of the natural projection $\pr_E:\Pi_E\twoheadrightarrow G_k$. 
We have a commutative diagram
$$
\CD
G_k  @>s>>   G_E \\
@V{\id}VV    @VVV  \\
G_k  @>{\tilde s}>> \Pi _E \\
\endCD
$$
where the right vertical map is the natural surjection.
The conjugacy class of the section $\tilde s$ corresponds to a unique element of $H^1(G_k,\Pi_{\overline E})$, which we will denote also $\tilde s\in
H^1(G_k,\Pi_{\overline E})$.
We have a natural exact sequence arising from Kummer theory
$$1\to E(k)^{\wedge}\to H^1(G_k,\Pi_{\overline E})\to TH^1(G_k,E)\to 1,\tag {$3$}$$
where $E(k)^{\wedge}\defeq \underset{n\ge 1}\to{\varprojlim} \frac {E(k)}{nE(k)}$ is the profinite completion of the finitely generated discrete group $E(k)$, 
$H^1(G_k,\Pi_{\overline E})$ is the profinite Galois cohomology group of the continuous $G_k$-module $\Pi_{\overline \Pi}$, 
and $TH^1(G_k,E)$ is the Tate module of the Galois cohomology group 
$H^1(G_k,E(\overline k))$.

\proclaim {Lemma 4.4} We use the same notations and assumptions as in 4.3.
The element $\tilde s\in H^1(G_k,\Pi_{\overline E})$, corresponding to the section $\tilde s:G_k\to \Pi_E$,
lies in the subgroup $E(k)^{\wedge}\subset H^1(G_k,\Pi_{\overline E})$.
\endproclaim

\demo {Proof} The existence of the section $s$ gives rise {\it naturally} to an adelic point $(x_v)_v\in E(\Bbb A_k)$ (cf. Proof of Proposition 1.4), where $x_v$ is uniquely 
determined at the finite places $v$. At a (possible) real place $v$ of $k$ only the connected component of $E(k_v)$ containing $x_v$ is well defined.
We have a natural commutative diagram of exact sequences
$$
\CD
1 @>>> E(k)^{\wedge} @>>>  H^1(G_k,\Pi_{\overline E}) @>>> TH^1(G_k,E) @>>> 1 \\
@.        @VVV              @VVV           @VVV    \\
1@>>> \prod _v E(k_v)^{\wedge}@>>> \prod _v H^1(G_{k_v},\Pi_{\overline E})@>>> \prod _v TH^1(G_{k_v},E) @>>> 1 \\
\endCD
$$
where the product in the lower exact sequence is taken over all places $v$ of $k$. The image of $\tilde s$ in $\prod _v H^1(G_{k_v},\Pi_{\overline E})$
is $(x_v)_v\in \prod _v E(k_v)^{\wedge}$.  The kernel of the natural map $TH^1(G_k,E)\to  \prod _v TH^1(G_k,E)$ is the Tate module
of the Shafarevich-Tate group of $E$; it is trivial if we assume that the latter is finite.  
Hence the image of $\tilde s$ in $TH^1(G_k,E)$ is trivial and $\tilde s$ lies in $E(k)^{\wedge}$ as claimed.
Note that in the above commutative diagram it is well-known that the left and middle vertical maps are injective.

Alternatively, the adelic point $(x_v)_v\in E(\Bbb A_k)$ survives every finite \'etale abelian descent obstruction
(cf. [Stoll], Definition 5.2), as follows easily from the existence of the global section $s$ (see also [Harari-Stix], Proposition 1.1), i.e. $(x_v)_v\in E(\Bbb A_k)^{\f-\ab}$
in the terminology of Stoll, where $E(\Bbb A_k)^{\f-\ab}$ is the set of adelic points cut out by the finite \'etale abelian conditions (cf. loc. cit. Definition 5.4).
This implies that $\tilde s$ lies in the Selmer group $\Sel (k,E)^{\wedge}\subset   H^1(G_k,\Pi_{\overline E})$ by a result of Stoll (cf. [Stoll], the discussion preceding 
Corollary 6.2). Furthermore, we have a natural identification $E(k)^{\wedge}\isom \Sel (k,E)^{\wedge}$, since we assumed the Shafarevich-Tate group of $E$ to be finite 
(cf. loc. cit.). Hence $\tilde s\in E(k)^{\wedge}$.
\qed
\enddemo

\definition {Remark 4.5} In fact, one can show slightly more than the statement in Lemma 4.4. For a finite closed subscheme $S\subset E$ denote by $J_S$ the corresponding generalised jacobian with modulus $S$. Write $J_S(k)^{\wedge}$ for the profinite completion of the group of $k$-rational points $J_S(k)$ of $J_S$. We have a natural homomorphism
$\varphi:\underset{S}\to{\varprojlim}\  J_S(k)^{\wedge}\to E(k)^{\wedge}$, where $\underset{S}\to{\varprojlim}\  J_S(k)^{\wedge}$ is the projective limit of the $J_S(k)^{\wedge}$'s. One can show that the element $\tilde s$ as in Lemma 4.4 lies in the image $\Im \varphi$ of the above homomorphism $\varphi$.
\enddefinition

In the framework of the above discussion one can characterise, using a result of Stoll, those sections $s:G_k\to G_E$ which arise from rational points as follows.

\proclaim {Proposition 4.6} We use the same notations and assumptions as in 4.3.
The image of the section $s:G_k\to G_E$ is contained in the decomposition group $D_x$ associated to a rational point
$x\in E(k)$ if and only if the induced section $\tilde s:G_k\to \Pi_E$ corresponds to an element $\tilde s\in H^1(G_k,\Pi_{\overline E})$ which lies in the subgroup $E(k)$ of $E(k)^{\wedge}\subset H^1(G_k,\Pi_{\overline E})$. (We already know that the element $s$ lies in $E(k)^{\wedge}$ by Lemma 4.4.)
\endproclaim

\demo {Proof} Note that the discrete group $E(k)$ naturally embeds into its profinite completion $E(k)^{\wedge}$, since it is finitely generated.
First, one easily observes that if the section $s$ arises from a rational point, i.e. if its image $s(G_k)\subset D_x$ is contained in a decomposition group 
associated to a rational point $x\in E(k)$, then the corresponding element $\tilde s \in H^1(G_k, \Pi_{\overline E})$ equals $s_x\in E(k)$, where $s_x\in H^1(G_k,\Pi_{\overline E})$ is the element corresponding to $x\in E(k)$. Here $E(k)$ is viewed as a subgroup of $H^1(G_k,\Pi_{\overline E})$ via the Kummer sequence, the natural map $E(k)\to E(k)^{\wedge}$ being injective.

Second, assume that $\tilde s=s_x\in E(k)$ for some rational point (necessarily unique) $x\in E(k)$. We will show that the image $s(G_k)\subset G_E$
of the section $s$ is contained in a decomposition group $D_x$ associated to the rational point $x$.
We use the following well-known argument in anabelian geometry.
In order to show that $s(G_k)\subseteq D_x$ it suffices to show (using a limit argument, and Faltings theorem on the finiteness of the set of rational points of a smooth, hyperbolic,
connected, and proper curve over a number field) that for every open subgroup $H$ of $G_E$ corresponding to a finite (possibly ramified) morphism $Y\to E$, where $Y$ has genus at least $2$, and such that $s(G_k)\subset H$,
we have $Y(k)\neq \varnothing$. Indeed, in this case the projective limit $\underset{Y}\to{\varprojlim} Y(k)$, where the limit is taken over all such $Y$'s, is non empty.
Consider a pro-point in $\underset{Y}\to{\varprojlim} Y(k)$, and its image $x'\in E(k)$. Then $s(G_k)\subseteq D_{x'}$, where $D_{x'}$ is a decomposition subgroup 
associated to $x'$. Moreover, $x'=x$ necessarily.

Next, let $H\subseteq G_X$ be an open subgroup corresponding to a finite morphism $g:Y\to E$, where $Y$ has genus at least $2$, and such that
$s(G_k)\subset H$. Then $H$ is naturally identified with the absolute Galois group $G_Y\defeq \Gal (K_E^{\sep}/K_Y)$, and the section $s$ restricts to a section
$s:G_k\to G_Y$ of the natural projection $G_Y\twoheadrightarrow G_k$. Similar arguments as the one used in the proof of Proposition 4.2 
imply that the existence of the section $s:G_k\to G_Y$ gives rise to an adelic point $(y_v)_v \in Y(\Bbb A_k)$ 
which survives every finite \'etale abelian descent obstruction, i.e. $(y_v)_v\in Y(\Bbb A_k)^{\f-\ab}\neq \varnothing$. Moreover, the image of $(y_v)_v$ in $E(\Bbb A_k)$
via the natural map $Y(\Bbb A_k)\to E(\Bbb A_k)$, which is 
induced by the natural morphism $g:Y\to E$, coincides with the adelic point $(x_v)_v\in E(\Bbb A_k)$ arising from the rational point $x\in E(k)$.
Let $Z\subset Y$ be the preimage (as a subscheme) of the rational point $x\in E(k)$. Thus, $Z$ is a finite $k$-scheme, and $(y_v)_v \in Z(\Bbb A_k)$.
We have $Z(k)=Z(\Bbb A_k)\cap Y(\Bbb A_k)^{\f-\ab}$ by a result of Stoll (cf. [Stoll], Theorem 8.2). In particular, $Z(k)\neq \varnothing$. Thus, $Y(k)\neq \varnothing$.
 This finishes the proof of Proposition 4.6.
\qed
\enddemo

In fact, the validity of the BGASC for elliptic curves over number fields gives a characterisation of the discrete group of rational points of an elliptic curve inside its profinite completion.
More precisely, we have the following which follows easily from Proposition 4.6.

\proclaim {Proposition 4.7} We use the same notations as above. Let $E$ be an elliptic curve over a number field $k$. 
Assume that the BGASC holds true for $E$ (cf. Definition 1.1).
Let $\tilde s\in E(k)^{\wedge}$, which we view as an element of 
$H^1(G_k,\Pi_{\overline E})$. Then $\tilde s$ lies in $E(k)$ if and only if a corresponding section $\tilde s:G_k\to \Pi_E$ of the natural projection  $\Pi_E\twoheadrightarrow G_k$
can be lifted to a section $s:G_k\to G_E$ of the natural projection $G_E\twoheadrightarrow G_k$, i.e. if there exists a section $s:G_k\to G_E$ and a commutative diagram
$$
\CD
G_k @>s>>  G_E \\
@V{\id}VV      @VVV    \\
G_k @>\tilde s>>  \Pi_E \\
\endCD
$$
where the right vertical map is the natural surjection.
\endproclaim

Proposition 4.7 implies immediately the following which was observed by Stoll (cf. [Stoll], Remark 8.9). See also [Harari-Stix] Theorem 3.5.

\proclaim {Corollary 4.8} Let $k$ be a number field, and $E$ an elliptic curve over $k$. Assume that the Shafarevich-Tate group
of $E$ is finite, and that the rank of the Mordell-Weil group of $E$ is trivial, i.e. $E(k)$ finite. Then the BGASC holds true for $E$.
\endproclaim

\subhead
\S 5. Birational Sections for Genus $g\ge 2$ Curves over Number Fields
\endsubhead

In this section we will establish some observations on the BGASC in the case of genus $g\ge 2$ curves over number fields.

\definition {5.1} Assume that $X$ is a proper, smooth, {\bf hyperbolic}, and geometrically connected {\bf curve} over a {\bf number field} $k$.
Assume that the {\bf Shafarevich-Tate group} of the jacobian $J\defeq J_X$
of $X$ is {\bf finite}. 
Let $s:G_k\to G_X$ be a section of the natural projection $G_X\twoheadrightarrow G_k$. Then $X$ has a rational divisor of degree $1$ (cf. [Esnault-Wittenberg]), and we can embed $X$ into $J$. Let $\Pi _{J}$ be the arithmetic fundamental group of $J$ which sits naturally in an exact sequence 
$$0\to T \overline J\to \Pi_J \to G_k\to 1,$$
where $T \overline J$ is the Tate module of $\overline J\defeq J\times _k \bar k$. Thus, $\Pi_J$ corresponds naturally to the quotient of $G_X$ which is the geometrically abelian \'etale fundamental group of $X$.
We fix a base point of the torsor of splittings of the above exact sequence which arises from the splitting associated to the zero section.
Recall the Kummer exact sequence 
$$0\to J(k)^{\wedge}\to H^1(G_k,T\overline J)\to TH^1(G_k,J)\to 1.$$
Similar arguments used in the proof of  Proposition 4.6 yield the following. 
\enddefinition

\proclaim {Proposition 5.2} We use the same notations and hypothesis as in 5.1. Let $s:G_k\to G_X$ be a section of the natural projection $G_X\twoheadrightarrow G_k$,
$\tilde s:G_k\to \Pi_J$ the section of the natural projection $\Pi_J\twoheadrightarrow G_k$ which is naturally induced by $s$, and $\tilde s\in H^1(G_k,T\overline J)$ the corresponding element of $H^1(G_k,T\overline J)$. Then $\tilde s$ lies in the subgroup $J(k)^{\wedge}$ of $H^1(G_k,T\overline J)$. Moreover, the image $s(G_k)\subset G_X$ of the section $s$ is contained in the decomposition group $D_x$ associated to a rational point
$x\in X(k)$ if and only if the above elements $\tilde s\in J(k)^{\wedge}$ lies in the subgroup $J(k)$ of $J(k)^{\wedge}$. 
\endproclaim

One can deduce, as a consequence of Proposition 5.2, the following.

\proclaim {Proposition 5.3} We use the same notations and hypothesis as in 5.1. Let $s:G_k\to G_X$ be a section of the natural projection $G_X\twoheadrightarrow G_k$,
and for each place $v$ of $k$ denote by $s_v:G_{k_v}\to G_{X_v}$ the corresponding section of the natural projection $G_{X_v}\twoheadrightarrow G_{k_v}$ (cf. proof of proposition 1.4).
Let $x\in X(k)$ be a rational point. Assume that for each place $v$ of $k$ the section $s_v$ arises from $x\in X(k)\subset X(k_v)$. In other words the image $s_v(G_{k_v})\subset 
\Tilde D_{x}$ is contained in a decomposition group $\Tilde D_{x}\subset G_{X_v}$ associated to the point $x \in X(k_v)$.
Then the section $s$ arises from the rational point $x$, i.e. the image $s(G_k)\subset D_{x}$ is contained in a decomposition group $D_{x}\subset G_X$ associated to the rational point $x$.
\endproclaim

\demo{Proof}
Indeed, with the same notations as in Proposition 5.2, in this case we have $\tilde s=x$ as an element of $X(k)\subset J(k)\subset J(k)^{\wedge}$.
\qed
\enddemo

\definition {Remark 5.4} The above discussion in the case of a curve $X$ of genus at least $2$ is related to the adelic-Mordell conjecture of Stoll
(cf. [Stoll]), which predicts that inside $\prod_vJ(k_v)$ the intersection $J(k)^{\wedge}\cap \prod _v X(k_v)$ is exactly $X(k)$. In fact the validity of Stoll's conjecture would imply, with the notation in Proposition 5.2, that $\tilde s$ lies automatically in $J(k)$, hence the validity of the BGASC for $X$ would follow.
However, in the case of an elliptic curve, Proposition 4.6 does not seem to be a priori related to Stoll's conjecture and the results in [Stoll]. 
\enddefinition

Finally, we observe the following.

\proclaim {Lemma 5.5} Let $k$ be a number field and $X$ a proper, smooth, and geometrically connected curve over $k$. Assume that there exists
an elliptic curve $E$ over $k$ with trivial Shafarevich-Tate group and with trivial Mordell-Weil rank. Then there exists a finite morphism
$f:X'\to X$ of degree $\deg (f)\le 2$ such that the BGASC holds true for $X'$ as a $k$-curve.
\endproclaim

\demo {Proof} Let $E$ be an elliptic curve over $k$ with trivial Shafarevich-Tate group and with trivial Mordell-Weil rank. Let $\tilde g:X\to \Bbb P^1_k$
be a finite morphism, and $\tilde f:E\to \Bbb P^1_k$ a morphism of degree $2$. Let $X'\defeq X\times _{\Bbb P^1_k} E$. We have a commutative diagram
$$
\CD
X' @>g>>  E \\
@V{f}VV      @V{\tilde f}VV    \\
X @>{\tilde g}>>  \Bbb P^1_k \\
\endCD
$$
where $f:X'\to X$ is a finite morphism of degree $\deg (f)\le 2$.
We choose the function $\tilde g$ so that $X'$ is geometrically connected.
We have a commutative diagram of exact sequences of absolute Galois groups
$$
\CD
1 @>>> G_{\overline X'} @>>>  G_{X'}@>\pr_X>> G_k  @>>>  1 \\
@.        @VVV              @VVV           @V{\id}VV    \\
1@>>> G_{\overline E}@>>> G_{E} @>{\pr_{\Bbb P^1_k}} >> G_k  @>>>  1\\
\endCD
$$
where the left and vertical maps are natural inclusions. Let $s:G_k\to G_{X'}$ be a group-theoretic section of the natural projection
$G_{X'}\twoheadrightarrow G_k$. Then $s$ induces naturally a group-theoretic section  $s':G_k\to G_E$ of the natural projection
$G_E\twoheadrightarrow G_k$. Moreover, one observes easily that the section $s$ arises from a rational point $x'\in X'(k)$ if and only if the section
$s'$ arises from a rational point $x\in E(k)$.  The section $s':G_k\to G_E$ arises from a rational point $x\in E(k)$ by Corollary 4.8.
\qed
\enddemo

\definition {Remark 5.6} Recently it was proven by Mazur and Rubin (cf. [Mazur-Rubin]) that over any number field $k$ there exist elliptic curves over $k$
with trivial Mordell-Weil rank. As a consequence, Lemma 5.5 implies that for any curve $X$ over $k$ there exists  
a  finite morphism $f:X'\to X$ of degree $\deg (f)\le 2$ such that the BGASC holds true for $X'$, under the assumption that the Shafarevich-Tate groups of elliptic curves over $k$ are finite.
\enddefinition

$$\text{References.}$$
\noindent
[Esnault-Wittenberg] Esnault, H., Wittenberg, O., On abelian birational sections, Journal of the American Mathematical society, 
Volume 23, Number 3, July 2010, Pages 713-724.

\noindent
[Grothendieck]  Grothendieck, A., Brief an G. Faltings, (German), with an english translation on pp. 285-293.
London Math. Soc. Lecture Note Ser., 242, Geometric Galois actions, 1, 49-58, Cambridge Univ. Press,
Cambridge, 1997.

\noindent
[Harari-Stix] Harari, D., and Stix J., Finite descent obstructions and fundamental exact sequence.  arXiv:1005.1302. 

\noindent
[Koenigsmann] Koenigsmann, J., On the section conjecture in anabelian geometry.  J. Reine Angew. Math.  588 
(2005), 221--235.

\noindent
[Manin] Manin, Y. I., le groupe de Brauer-Grothendieck en g\'eom\`etrie diophantienne, Actes du congr\`es international des math\'ematiciens
(Nice, 1970), Tome 1, pp. 401-411. Gauthier-Villars, Paris (1971).  

\noindent
[Mazur-Rubin]  Mazur, B., Rubin, K., Ranks of twists of elliptic curves and Hilbert's tenth problem. Invent. Math. 181 (2010), no. 3, 541-575.

\noindent
[Pop] Pop. F., On the birational $p$-adic section Conjecture, Compos. Math. 146 (2010), no. 3, 621-637.

\noindent
[Sa\"\i di] Sa\"\i di, M., Around the Grothendieck anabelian section conjecture. London Math. Soc. lecture Note Ser. 393,
Non-abelian Fundamental Groups and Iwasawa Theory, 72-106, Cambridge Univ. Press, Cambridge, 2011.
Edited by John Coates, Minhyong Kim, Florian Pop, Mohamed Sa\"\i di, and Peter Schneider.

\noindent
[Stoll] Stoll, M., Finite descent obstructions and rational points on curves. Algebra Number Theory 1 (2007), no. 4, 349-391.

\bigskip

\noindent
Mohamed Sa\"\i di

\noindent
College of Engineering, Mathematics, and Physical Sciences

\noindent
University of Exeter

\noindent
Harrison Building

\noindent
North Park Road

\noindent
EXETER EX4 4QF 

\noindent
United Kingdom

\noindent
M.Saidi\@exeter.ac.uk

\end
\enddocument